\documentclass[12pt]{amsart}
\usepackage{amsmath, amsthm, amsfonts, amssymb, xypic, float}
\usepackage[mathscr]{euscript}
\usepackage{tikz}

\usetikzlibrary{arrows}
\restylefloat{table}

\newtheorem{thm}{Theorem}[section]

\newtheorem{defn}[thm]{Definition}
\newtheorem{conj}[thm]{Conjecture}

\newcommand{\cc}{\mathcal{C}}
\newcommand{\cs}{\mathcal{S}}

\newcommand{\ca}{\mathcal{A}}

\newcommand{\cb}{\mathcal{B}}

\newcommand{\ck}{\mathscr{K}}
\newcommand{\ct}{\mathscr{T}}
\newcommand{\tk}{\tilde{K}}
\newcommand{\tf}{\tilde{F}}
\newcommand{\C}{\mathbb{C}}
\newcommand{\Z}{\mathbb{Z}}

\newcommand{\codiag}{\underrightarrow{\mathrm{Dia}}\mathrm{g}}
\newcommand{\diag}{\underleftarrow{\mathrm{Dia}}\mathrm{g}}
\newcommand{\colim}{\underrightarrow{\mathrm{lim}}}
\newcommand{\llim}{\underleftarrow{\mathrm{lim}}}

\def\Prim{\mathrm{Prim}}
\def\Spec{\mathrm{Spec}}
\def\KHaus{\mathbf{KHaus}}
\def\Ab{\mathbf{Ab}}

\def\Mdot{\text{ .}}

\title{From Topology to Noncommutative Geometry: $K$-theory}
\author{Nadish de Silva}

\address{
Quantum Group, Department of Computer Science, University of Oxford \\
Wolfson Building\\
Parks Road\\
Oxford OX1 3QD \\
UK.}
   \email{nadish.desilva@utoronto.ca}

\date{}

\begin{document}

\pagestyle{empty}
\pagenumbering{gobble}
\maketitle

\begin{abstract}
We associate to each unital $C^*$-algebra $\ca$ a geometric object---a diagram of topological spaces representing quotient spaces of the noncommutative space underlying $\ca$---meant to serve the role of a generalized Gel'fand spectrum.  After showing that any functor $F$ from compact Hausdorff spaces to a suitable target category $\cc$ can be applied directly to these geometric objects to automatically yield an extension $\tilde{F}$ which acts on all unital $C^*$-algebras, we compare a novel formulation of the operator $K_0$ functor to the extension $\tilde K$ of the topological $K$-functor.
\end{abstract}

\section{Introduction}

Gel'fand duality provides a link between geometry and algebra by asserting that every unital, commutative $C^*$-algebra is isomorphic to $C(X)$, the algebra of continuous, complex-valued functions on a compact Hausdorff topological space $X$.  Indeed, the duality provides an equivalence between these categories.  This in turn yields a method of importing concepts of topology into the study of $C^*$-algebras: topological concepts are translated into the language of commutative algebra and generalized to not-necessarily-commutative algebras \cite{khalkhali}.  As a simple example, the open sets of a space $X$ are in ordered bijection with the closed, two-sided ideals of $C(X)$ and so, such ideals can be considered the noncommutative analogue of open sets for not-necessarily-commutative algebras.  A more sophisticated example comes from vector bundles over a space $X$ which can be identified with the finite, projective modules over $C(X)$ \cite{swan}.  Pursuing this allows one to define a generalization of topological $K$-theory, which is formulated in terms of vector bundles \cite{atiyah}, to all $C^*$-algebras and obtain operator $K$-theory, an extremely potent invariant \cite{elliott, rordam2}.  Gel'fand duality justifies considering $C^*$-algebras to be generalizations of topological spaces; this is the starting point of the program of noncommutative operator geometry~\cite{connes}.  The motivating metaphor of noncommutative geometry is to view a $C^*$-algebra as being the algebra of functions on a topological space.  Although, in the case of commutative $C^*$-algebras, this space is explicitly given by the Gel'fand spectrum, in the noncommutative case, the existence of such an underlying space is understood to be merely a convenient fiction which justifies the algebraic generalization of geometric concepts. 

In this note, we present a geometric construction meant to serve as a generalization of the Gel'fand spectrum to noncommutative \mbox{$C^*$-algebras}.  We outline an automatic method for the generalization of topological concepts to the setting of $C^*$-algebras which is distinct from the one described above.  The method requires essentially two steps.  First, we contravariantly associate to each unital $C^*$-algebra $\ca$ the geometric object $G(\ca)$ which is meant to play the role of the hypothetical space $X$ on which $\ca$ is the algebra of continuous, complex-valued functions.  Next, we apply directly to $G(\ca)$ the topological functor which we aim to extend.  Thus, we begin with any functor $F$ from the category of compact Hausdorff spaces to a suitable target category $\cc$ and we construct a functor $\tilde{F}$ from the category of unital $C^*$-algebras to $\cc$ such that $\tilde{F}$, restricted to commutative algebras, is naturally isomorphic to $F$ composed with the Gel'fand spectrum functor.  The results of this method of extending topological concepts can be compared with the constructions achieved via the canonical translation process.  Agreement of these two methods would justify regarding $G(\ca)$ as the underlying space of a noncommutative algebra.

Akemann and Pedersen \cite{pedersen} proposed to replace the translation process by working directly with Giles-Kummer's \cite{giles} and Akemann's \cite{akemann} noncommutative generalizations of the basic topological notions of open and closed sets.  In contrast, we do not employ algebraic generalizations of basic topological notions, but rather, we work with objects which slightly generalize the notion of topological space and come readily equipped with an alternative to the translation process.

In addition to the work of Akemann-Pedersen and Giles-Kummer on noncommutative generalizations of Gel'fand duality, there have been a number of  alternative approaches by authors including Alfsen \cite{Alfsen}, Bichteler et al. \cite{Bichteler},  Dauns-Hofmann \cite{Dauns}, Fell \cite{Fell}, Kadison \cite{Kadison}, Kruml et al. \cite{Kruml}, Krusy\'nski-Woronowicz \cite{Krusynski}, Resende \cite{Resende}, Schultz~\cite{Schultz}, and Takesaki \cite{Takesaki}.  An excellent discussion of many of these works is contained in a paper by Fujimoto \cite{Fujimoto}.

We begin by constructing the geometric object $G(\ca)$ we associate to a unital $C^*$-algebra $\ca$, and giving a motivation for it as a means of keeping track of all the quotient spaces of the noncommutative space underlying $\ca$ which happen to be genuine topological spaces.  After showing how this association of geometric objects to algebras allows us to automatically extend topological functors, we then relate the resulting extension of topological $K$-theory to a novel formulation of the operator $K_0$ functor.  

We then report on progress on how the basic topological notion of lattice of closed sets might be extended to the notion of lattice of closed, two-sided ideals, thereby establishing a conjectured relationship between the topologies of the geometric object $G(\ca)$ and $\Prim(\ca)$, the primitive ideal space of $\ca$. (The primitive ideal space of $\ca$ is the set of ideals which are kernels of irreducible $*$-representations of $\ca$ equipped with the hull-kernel topology; in the commutative case, the primitive ideal space coincides with the Gel'fand spectrum.)  Establishing this conjecture would allow considering $G$ to be an enrichment of $\Prim$, a  $C^*$-algebraic variant of the ring-theoretic spectrum functor $\Spec$ which is the basis of the sheaf-theoretic study of rings.  This might open up the possibility of investigating the use of sheaf-theoretic methods in noncommutative geometry.

\section{Spatial diagrams}

To each unital $C^*$-algebra $\ca$, we will associate a diagram of topological spaces $G(\ca)$ which is expressed formally as a contravariant functor with, as its codomain, the category of compact Hausdorff spaces.  Not only will these spaces vary with $\ca$ but the shape of the diagram (the domain of the functor) will as well.  This necessitates introducing the following constructions:

\begin{defn}For any category $\cc$, $\codiag(\cc)$, the covariant category of all diagrams in $\cc$, has as objects all the functors $D$ from any small category $\cs$ to $\cc$.  Morphisms from $D_1: \cs_1 \rightarrow \cc$ to $D_2: \cs_2 \rightarrow \cc$ are given by pairs $(f, \eta)$ where $f$ is a functor from $\cs_1 \rightarrow \cs_2$ and $\eta$ is a natural transformation from $D_1$ to $D_2 \circ f$. 

The contravariant category of all diagrams $\diag(\cc)$ has all contravariant functors to $\cc$ as objects; the morphisms from $D_1$ to $D_2$ are pairs $(f, \eta)$ where $f$ is a functor from $\cs_2 \rightarrow \cs_1$ and $\eta$ is a natural transformation from $D_1 \circ f$ to $D_2$. 
\end{defn}

The composition $(g,\mu) \circ (f, \eta)$ of two $\codiag(\cc)$-morphisms is given by $(g  f, (\mu f)\eta)$ where $(\mu f)_a$ is $\mu_f(a)$ and the composition of natural transformations is componentwise.  Note that if $F$ is a functor from $\cc$ to $\cc'$, $F$ naturally induces a functor from $\codiag(\cc)$ to $\codiag(\cc')$ which we will also denote by $F$.  Explicitly, if $D: \ca \rightarrow \cc$, then $F(D)$ is simply $F \circ D$.  For a $\codiag(\cc)$-morphism $(f, \eta)$, $F$ sends $(f, \eta)$ to the $\codiag(\cc')$-morphism $(f, F \eta)$ where $(F \eta)_a$ is $F(\eta_a)$.  The functor $F$ also induces, in a similar fashion, a functor from $\diag(\cc)$ to $\diag(\cc')$.  If $F$ is contravariant, then it induces a contravariant functor from $\codiag(\cc)$ to $\diag(\cc')$ and one from $\diag(\cc)$ to $\codiag(\cc')$.

\subsection{The spatial diagram of an algebra}

We can now construct the spatial diagram $G(\ca)$ associated to a unital $C^*$-algebra $\ca$.  First, we associate to $\ca$ a subcategory $S(\ca)$ of unital, commutative $C^*$-algebras with unital $*$-homomorphisms.  The objects of $S(\ca)$ are the unital, commutative sub-$C^*$-algebras of $\ca$.  For every inner automorphism $\alpha$ of $\ca$, that is, one acting by conjugation by a unitary $u \in \ca$,  object $U \subset \ca$, and any object $V \subset \ca$ containing $\alpha(U)$, there is a restriction $\alpha|_U$ from $U$ to $V$; these are the morphisms of $S(\ca)$.  Denote by $g(\ca)$ the inclusion functor from $S(\ca)$ to $\mathbf{uComC^*}$.  The diagram $G(\ca)$, an object of $\diag(\KHaus)$, is the Gel'fand spectrum functor $\Sigma$ composed with $g(\ca)$.

The association of spatial diagrams to algebras is contravariantly functorial.  To a unital *-morphism $\phi: \ca \rightarrow \cb$, we define $g(\phi)$ as the $\codiag(\mathbf{uComC^*})$-morphism $(f, \eta)$.  The functor $f$ maps a commutative subalgebra $U \subset \ca$ to $\phi(U) \subset \cb$; the action of $f$ on the restriction of an automorphism described by a unitary $u$ is to send it to the restriction of an automorphism described by $\phi(u)$.  The component of $\eta$ associated to $U \subset \ca$ is $\phi|_U$.  Finally, the $\diag(\KHaus)$-morphism $G(\phi)$ is $\Sigma(g(\phi))$.

The topological spaces in the diagram $G(\ca)$ should be thought of as being those which arise as quotient spaces of the hypothetical noncommutative space underlying $\ca$.  To see this, note that a sub-$C^*$-algebra of $C(X)$ yields an inclusion $i: C(Y) \rightarrow C(X)$ which corresponds to a continuous surjection $\Sigma i: X \rightarrow Y$; this surjection is a quotient map since both $X$ and $Y$ are compact and Hausdorff.  Thus, in accordance with the central tenet of noncommutative geometry, unital sub-$C^*$-algebras of a unital noncommutative algebra $\ca$ are to be understood as having an underlying noncommutative space which is a quotient space of the noncommutative space underlying $\ca$.  By considering only the commutative subalgebras, we are restricting our attention to the tractable quotient spaces: those which are genuine topological spaces.  The morphisms of the diagram---particularly, the images under the Gel'fand spectrum functor of inclusions---serve to track how these spaces fit together inside the noncommutative space underlying $\ca$.

\section{Extensions}

For a unital $C^*$-algebra $\ca$, we have constructed the diagram of topological spaces $G(\ca)$.  Any functor $F: \KHaus \rightarrow \cc$, from the category of compact Hausdorff spaces to $\cc$, can be applied directly to $G(\ca)$ to yield $F \circ G (\ca)$, a diagram in $\cc$.  When $F$ is contravariant (resp. covariant) and $\cc$ is cocomplete (resp. complete), we can take the colimit (resp. limit) of this diagram to yield a single object of $\cc$.  Combining all these steps gives us a new functor $\tilde{F}$ which acts on all unital $C^*$-algebras.  Since, when $\ca$ is commutative, $\ca$ is a terminal object in the category $S(\ca)$, it can be shown that $\tilde{F}$, restricted to commutative algebras, is naturally isomorphic to $F$ composed with the Gel'fand spectrum functor.  Intuitively, one should think of the extension process as decomposing a noncommutative space into its quotient spaces, retaining those which are genuine topological spaces, applying the topological functor to each one, and pasting together the result.

\begin{defn} 
For a contravariant functor $F: \KHaus \rightarrow \cc$ with $\cc$ cocomplete, $\tilde{F}: \mathbf{uC^*} \rightarrow \cc$ the \emph{extension of $F$} is $\colim \: F  G$.
\end{defn}

This definition requires the generalized colimit functor \linebreak $\colim: \codiag(\cc) \rightarrow \cc$.  It assigns to a functor $F: \ca \rightarrow \cc$ the same object of $\cc$ which is assigned to $F$ by the canonical colimit functor $\colim: \cc^\ca \rightarrow \cc$.  If $\eta$ is a natural transformation between $F$ and $G: \ca \rightarrow \cc$ then $\colim$ assigns to the $\codiag(\cc)$-morphism between $F$ and $G$ given by $(id_\ca, \eta)$ the same $\cc$-morphism assigned to $\eta$ by the colimit functor of $\cc^\ca$.  What is novel is the ability to assign $\cc$-morphisms between colimits of diagrams of different shapes.  A generalized limit functor $\llim: \diag(\cc) \rightarrow \cc$ also exists and allows defining the extension of a covariant functor $F$ as $\llim \: F  G$.

In the sequel of this section, we begin by explicitly describing how the generalized $\colim$ and $\llim$ functors act on arbitrary morphisms and demonstrating the functoriality of these generalized notions of universal constructions.  We then give an explicit description of the $\colim$ functor on the category of all diagrams of abelian groups which we will require in our formulation of operator $K$-theory.  We conclude the section by showing that, given a functor $F: \KHaus \rightarrow \cc$, its extension $\tilde{F}$, once restricted to commutative algebras, is naturally isomorphic to $F$ composed with the Gel'fand spectrum functor $\Sigma$.
\subsection{The generalized limit and colimit functors}
 Recall that the colimit of a functor $F$ from $\ca$ to a cocomplete category $\cc$ can be expressed as a coequalizer of two coproducts \cite[p355]{maclane}:
$$
   \xymatrix@+=3pc{\coprod_{u:i \rightarrow j}{F(\emph{\emph{dom}}u)} \ar@<0.5ex>[r]^-\theta \ar@<-0.5ex>[r]_-\tau & \coprod_{i}{F(i)}
   }
 $$
The first coproduct is over all arrows $u: i \rightarrow j$ of $\ca$ and the second is over all objects of $\ca$.  We denote the canonical injections for these coproducts by $$\lambda_u: F(\mathrm{dom}u) \rightarrow \coprod_{u:i \rightarrow j}{F(\mathrm{dom}u)}$$ $$\kappa_i: F(i) \rightarrow \coprod_{i}{F(i)} \Mdot$$  The morphisms $\theta$ and $\tau$ can be defined by specifying their compositions with the $\lambda_u$:
$$
\theta \lambda_u = \kappa_{\emph{\emph{dom}}u} \>\>\>\>\>\>\>\>\>\>\>\>\>\>\>
\tau \lambda_u = \kappa_{\emph{\emph{cod}}u} F(u)
$$
The advantage of this coequalizer presentation of the colimit of is that we may determine a $\cc$-morphism between the colimits of two functors $F$ and $G$ by specifying a natural transformation between their coequalizer diagrams.  That is, by giving its components;  $\cc$-morphisms $N$ and $M$ such that the following diagrams commute:
$$
   \xymatrix{\coprod_{u:i \rightarrow j}{F(\emph{\emph{dom}}u)} \ar@<0.5ex>[r]^-\theta  \ar[d]^N & \coprod_{i}{F(i)} \ar[d]^M \\
\coprod_{u':i' \rightarrow j'}{G(\emph{\emph{dom}}u')} \ar@<0.5ex>[r]^-{\theta'}  & \coprod_{i'}{G(i')}}
$$
 
$$
   \xymatrix{\coprod_{u:i \rightarrow j}{F(\emph{\emph{dom}}u)}  \ar@<-0.5ex>[r]_-\tau \ar[d]^N & \coprod_{i}{F(i)} \ar[d]^M \\
\coprod_{u':i' \rightarrow j'}{G(\emph{\emph{dom}}u')} \ar@<-0.5ex>[r]_-{\tau'} & \coprod_{i'}{G(i')}}
    $$
We denote the canonical injections into the coproducts for $G$ by $\lambda_u '$ and $\kappa_i '$.

Given a $\codiag(\cc)$-morphism $(f, \eta)$ between $F$ and $G$ we define $N$ and $M$ by giving their compositions with the canonical injections:
$$
N \lambda_u = \lambda_{f(u)}' \eta_{\emph{\emph{dom}}u} \>\>\>\>\>\>\>\>\>\>\>\>\>\>\>
M \kappa_i = \kappa_{f(i)}' \eta_i
$$
It is straightforward to verify that $\theta ' N = M \theta$ and that $\tau ' N = M \tau$ by computing the composition of these maps with the $\lambda_u$.  The \linebreak $\cc$-morphism assigned by $\colim$ to $(f, \eta)$ is then defined to be that morphism which is induced by the natural transformation between the coequalizer diagrams for the colimits of $F$ and $G$ whose components are $N$ and $M$ .  Functoriality of $\colim$ is then straightforwardly verified by computing the compositions of the components of the natural transformations induced by $(f, \eta)$ and $(g, \mu)$ and seeing that the resulting natural transformation is the same as the one induced by $(gf, (\mu f)\eta)$.
  
A dual construction exists in the case where $\cc$ is complete, in which case we have a generalized limit functor $\llim: \diag(\cc) \rightarrow \cc$.

The generalized colimit construction is best illustrated by the example of abelian groups.  Let $F: A \rightarrow \Ab$ and $G: B \rightarrow \Ab$ be two diagrams in $\codiag(\Ab)$ and $(f, \eta)$ be a morphism from $F$ to $G$.  

First, we describe the colimit of $F$ in $\Ab$ (and thus in $\codiag(\Ab)$).  Let $d$ be the direct sum of the groups $F(a)$ over all objects $a$ in $A$.  If $g$ is an element of the $F(a)$, we use the notation $(g)_a$ to indicate the element of $d$ which is $g$ in the $a^{th}$ component and 0 in all others.  The colimit of $F$ is $d$ modulo an equivalence relation: for any arrow $u: a_1 \rightarrow a_2$ of $A$, $(g)_{a_1}$ is identified with $(F(u)(g))_{a_2}$.

To define how the functor $\colim: Diag(\Ab) \rightarrow \Ab$ acts on $(f, \eta)$, it is enough to say how the group homomorphism $\colim((f, \eta))$ acts on elements of the colimit of $F$ of the form $[(g)_a]$.  The image of such an element under $\colim((f, \eta))$ is $[(\eta_a(g))_{f(a)}]$.  This is well defined for if $u: a_1 \rightarrow a_2$ identifies $(g)_{a_1}$ with $(F(u)(g))_{a_2}$, then $G \circ f(u)$ identifies $(\eta_{a_1}(g))_{f(a_1)}$ with $(\eta_{a_2}(F(u)(g)))_{f(a_2)}$.

\subsection{Extensions and the commutative case}

When  $\ca$ is commutative, $\ca$ is terminal in the category $S(\ca)$ as there is an inclusion morphism from every unital, commutative subalgebra of $\ca$ into $\ca$.  As a consequence, the diagram $F(G(\ca))$ has  $F(\Sigma(\ca))$ as a terminal object when $F$ is contravariant.  Since, for any diagram whose domain has a terminal object, the canonical injection from the terminal object to the colimit is an isomorphism, we have that $\tf(\ca) \simeq F \circ \Sigma(\ca)$.  For a unital $*$-homomorphism $\phi: \ca \rightarrow \cb$ between commutative algebras, $\tf(\phi)$, the isomorphisms , and $F \circ \Sigma(\phi)$ together form a commutative square.  Thus, these isomorphisms define a natural equivalence between $\tf|_{ComC^*}$ and $F \circ \Sigma$.  A dual argument holds in the case where $F$ is covariant.

\section{$K$-theory, Topological to $C^*$-algebraic}

Given a compact Hausdorff space $X$, the isomorphism classes of vector bundles over $X$ form an abelian monoid $V(X)$ under direct sum.  The universal abelian group of $V(X)$ is the $K$-group of $X$ and is denoted by $K(X)$.  The topological $K$-functor is the basis of the extraordinary cohomology theory $K$-theory which is of considerable importance in geometry.  It is generalized to arbitrary $C^*$-algebras using the canonical translation process.  As vector bundles over $X$ are, by the Serre-Swan theorem, in correspondence with finite, projective \linebreak $C(X)$-modules, we define $V(\ca)$ to be the monoid of isomorphism classes of finite, projective $\ca$-modules with direct sum and find that it, too, is abelian.  Again using the universal group construction yields the operator $K_0$ functor.  We aim to define operator $K_0$ directly in terms of the topological $K$-theory of spatial diagrams.

The first step in computing the $K_0$-group of a $C^*$-algebra is to take its stabilization.  The stabilization functor $\ck: \mathbf{C^*Alg} \rightarrow \mathbf{C^*Alg}$ acts on an algebra $\ca$ by taking it to $\ca \otimes \ck$, its tensor product with the \linebreak $C^*$-algebra $\ck$ of compact operators on a separable Hilbert space, and sends *-morphisms $\phi$ to $\phi \otimes id_{\ck}$.  This operation is idempotent: \linebreak $\ck \circ \ck \simeq \ck$.  The $K$-theory of a $C^*$-algebra $\ca$ is the same as its stabilization $\ck(\ca)$; in fact, $K_0 \circ \ck \simeq K_0$.

Although it remains open whether $K_0 \simeq \tk$ we find that our aim can be achieved by applying $K$ not to the diagrams $G$ but to subdiagrams $G_f$ after stabilizing.

\begin{thm}
$K_0 \circ \ck \simeq K_0 \simeq \tk_{f} \circ \ck$ as functors from unital \linebreak $C^*$-algebras to abelian groups.
\end{thm}

Here, $\tk_f$ is defined for unital $C^*$-algebras in a manner most similar to how $\tk$ is: as $\colim \: K \circ G_f$.  The only difference in the definitions of $G$ and $G_f$ is that the sub-$C^*$-algebras of $\ca$ which are the objects in the full subcategory $S_f(\ca)$ of $S(\ca)$ are, in addition to being unital and commutative, finite dimensional.  As all $C^*$-algebras isomorphic to $\ca \otimes \ck$ lack a unit, we must extend $\tk_f$ to all $C^*$-algebras and we do so using the same method used for $K_0$: if $\pi: \ca^+ \rightarrow \C$ is the projection from the unitalization of $\ca$ to the subalgebra of scalar multiples of the unit, then $\tk_f(\ca)$ is defined as the kernel of $\tk_f(\pi)$ \cite{fillmore}.

For a unital $C$*-algebra $\ca$, the $K_0$-group can be expressed in a rather simple fashion \cite{rordam}.  First, we define an equivalence relation on projections in $\ca \otimes \ck$: $P \sim Q$ if and only if there is a unitary $u \in \mathcal{(A \otimes \ck)}^+$ such that $P = uQu^*$.  

If $P$ and $P'$ are two  projections in $\ca \otimes \ck$, there exist projections $Q$ and $Q'$ such that $P \sim Q$, $P' \sim Q'$, and $Q$ is orthogonal to $Q'$.  The abelian addition operation on projections is defined by \linebreak $[P] + [P'] = [Q + Q']$.  We have thus defined the abelian monoid $V(\ca)$.

Then $K_0(\ca)$ is simply the collection of formal differences $[P] - [Q]$ of elements of $V(\ca)$ with $[P] - [Q] = [P'] - [Q']$ precisely when there exists $[R]$ such that $[P] + [Q'] + [R] = [P'] + [Q] + [R]$.  A morphism $\phi: \ca \rightarrow \cb$ between unital $C$*-algebras naturally induces a morphism $\ck(\phi)$ between their stabilizations.  When restricted to projections, this map is is well-defined and additive on the equivalence classes, yielding a morphism of abelian monoids $V(\phi): V(\ca) \rightarrow V(\cb)$.  By composing with the universal abelian group functor, we turn this abelian monoid homomorphism into an abelian group homomorphism between the $K_0$ groups of $\ca$ and $\cb$.

\begin{proof}
For unital $\ca$, we will compute $\tilde{K}_f \circ \ck(\ca) =  \colim \: K \circ G_f(\ca \otimes \ck)$ in stages and present it in such a form which makes it clear that it is isomorphic in a natural way to $K_0(\ca)$.  We begin by describing how the functor $\tilde{K}_f$ acts on unital algebras before generalizing to the non-unital case.

The objects of $S_f(\ca)$ are the finite-dimensional, unital, commutative sub-$C^*$-algebras of $\ca$.  The morphisms are given by the restrictions of inner automorphisms.  These morphisms are all of the form $i \circ r$ where $i$ is an inclusion and $r$ is an isomorphism between subalgebras which are related by unitary rotation.  

Under the Gel'fand spectrum functor, the image of such an object is a finite, discrete space whose points are in correspondence with the atomic projections of the subalgebra.  The images of the inclusions $i: U \rightarrow V $ are functions $\Sigma(i): \Sigma(V) \rightarrow \Sigma(U)$ with the property that whenever a point $p \in \Sigma(U)$ corresponds to a projection $P$ atomic in $U$, then $P$ is the sum of the atomic projections in $V$ which correspond to the points of $[\Sigma(i)]^{-1}(p)$.  The isomorphisms are sent by the spectrum functor to bijections which connect points corresponding to unitarily equivalent projections.  

Under the topological $K$-functor, each object yields a direct sum of copies of $\Z$: one for each point (a trivial vector bundle of each dimension and formal inverses).  Taking the colimit of the diagram then yields, as described in the previous section, a direct sum of the groups indexed by the objects of $S_f(\ca)$ modulo the relations generated by the morphisms of the diagram.  In our case, this is a quotient of the direct sum of a copy of $\Z$ for each pair $(U, P)$ where $U$ is a finite-dimensional, unital, commutative sub-$C^*$-algebra of $\ca$ and $P$ is an atomic projection in $U$.

The image under $K \circ \Sigma$ of the inclusions result in the identifications        of a generator of a copy of $\Z$ associated to a pair $(U, P)$ with the sums of generators associated to pairs $(V, P_i)$ whenever $U \subset V$ and $P_i$ sums to $P$.  Every projection $P \in \ca$ is an atomic projection in the subalgebra $\C P + \C (1 - P)$  which is included in every subalgebra which contains $P$ as a member.  As the generators associated to the same projection $P$ atomic in different subalgebras are all identified in the colimit, we see that we may speak of the element of the colimit group $[(P)]$ associated to $P$ without reference to which subalgebra it appears in.  Thus, the abelian group $\tilde{K}_f(\ca)$ can be viewed as a quotient of the free abelian group generated by the elements $[(P)]$.  The first class of identifications consists of those between elements $[(P)]$ with the sum of elements $[(P_i)]$ whenever $P_i$ are mutually orthogonal and sum to $P$.  The isomorphisms in the diagram provide the second class; they ensure that the elements associated to unitarily equivalent projections are also identified.

For non-unital algebras such as $\ca \otimes \ck$, we need to determine the kernel of $\tilde{K}_f(\pi)$ with $\pi$ the canonical projection from $(\ca \otimes \ck)^+$ to $\C$.  This is not so difficult, however, as all projections in $(\ca \otimes \ck)^+$ are of the form $P$ or $1 - P$ for $P \in \ca \otimes \ck$.  As $[(1-P)] = [(1)] - [(P)]$ we see that all elements of the colimit group can be expressed using only elements associated either to the identity projection or to projections in $\ca \otimes \ck$.  It is precisely those elements of $\tilde{K}_f((\ca \otimes \ck)^+)$ which can be expressed using only elements associated to projections in $\ca \otimes \ck$ which are in the kernel.

We are ready to define the natural transformation $\eta: K_0 \rightarrow \tilde{K_f} \circ \ck$.  The component of the natural isomorphism $\eta_\ca$ sends $[P] - [Q]$ in $K_0(\ca)$ to $[(P)] - [(Q)]$.  This is easily seen to be well defined, for if \linebreak $[P] + [Q'] + [R] = [P'] + [Q] + [R]$ then we may find mutually orthogonal representatives of all these projections and demonstrate that \linebreak $[(P)] - [(Q)] = [(P')] - [(Q')]$.  Preservation of addition follows by a similar argument.

We define an inverse map to demonstrate bijectivity.  The element $[(P)] \in \tilde{K_f}(\ca \otimes \ck)$ is sent by $\eta_\ca^{-1}$ to $[P]$.  Since the equivalences induced by the morphisms of the diagram are respected at the level of $K$-theory, this map is well-defined.

To demonstrate the naturality of these isomorphisms, let \linebreak $\phi: \ca \rightarrow \cb$ be a unital *-morphism.  We will require the naturally induced *-morphism $\ck\phi: \ca \otimes \ck \rightarrow \cb \otimes \ck$ which is defined explicitly by $\phi \otimes id_\ck$. 
 \begin{equation}
   \xymatrix{{K_0(A)}  \ar@<-0.5ex>[r]_-{\eta_A} \ar[d]^{K_0(\phi)} & {\colim \tilde{K}_f \circ \ck(A)} \ar[d]^{\colim \tilde{K}_f \circ \ck(\phi)} \\
{K_0(B)} \ar@<-0.5ex>[r]_-{\eta_B} & {\colim \tilde{K}_f \circ \ck(B)}
   }
 \end{equation}
  
An element $[P] - [Q]$ in $K_0(\ca)$ is mapped to $[(P)] - [(Q)]$ by $\eta_\ca$.  This is then mapped to $[(\ck{\phi}(P))] - [(\ck{\phi}(Q))]$ by $\colim \: \tilde{K}_f \circ \ck(\phi)$.  Alternatively, $[P] - [Q]$ in $K_0(\ca)$ is mapped to $[\ck{\phi}(P)] - [\ck{\phi}(Q)]$ in $K_0(\cb)$ which is in turn mapped to $[(\ck{\phi}(P))] - [(\ck{\phi}(Q))]$ by $\eta_\cb$.

\end{proof}

\section{Conclusions and Future Work}

As we have seen, operator $K_0$ theory arises as the finite-dimensional fragment of the extension $\tk$ of topological $K$-theory.  Although this applies only to the stabilizations of unital $C^*$-algebras, this is no great loss; any algebra is stabilized before its $K_0$ group is computed.  It is also possible that the relationship between $K$ and $\tilde{K}$ is even closer than what we have proven: perhaps $K_0 \simeq \tilde{K}$.

In order to justify thinking of $G: uC^* \rightarrow \diag(\KHaus)$ as a noncommutative Gel'fand spectrum functor, we wish to find more examples of topological concepts which can be automatically translated to their noncommutative version via the extension method based on $G$.  We have considered the example of the notion of closed sets.  Let $\ct: \KHaus \rightarrow \mathbf{CMSLat}$ be the functor which assigns to a compact Hausdorff topological space its lattice of closed sets with containment and, to a continuous function, the homomorphism of complete meet-semilattices which sends a closed set to its image under the continuous function. We conjecture that:

\begin{conj} \label{cstarconj}
The extension of the  lattice of closed sets functor,  $\tilde{\ct}: \mathbf{uC}^* \rightarrow \mathbf{CMSLat}$, assigns to a unital $C^*$-algebra its lattice of closed two-sided ideals.
\end{conj}

  Equivalently, the hull-kernel topology on the $C^*$-algebra's primary ideal space can be recovered as the limit of the topologies of $G$.  The functor $\Prim$ is the $C^*$-algebraic variant of the spectrum functor $\Spec$ which assigns  to a commutative ring a space whose hull-kernel topology provides the basis for sheaf-theoretic ring theory. 

For the purposes of proving the conjecture, we find it useful to rephrase the problem.

\begin{defn}
A \emph{partial ideal}  \cite{reyes} of a unital $C^*$-algebra $\ca$ is a choice of closed ideal $I_V$ from each unital, commutative sub-$C^*$-algebra $V \subset \ca$ such that whenever $V \subset V'$, the ideal $I_V$ can be recovered from $I_{V'}$ as $I_V = I_{V'} \cap V$. 
\end{defn}

Every closed, two-sided ideal of $\ca$, which we will call a \emph{total ideal}, gives rise to a partial ideal in a natural way by choosing $I_V$ to be $I \cap V$.  Elements of the lattice $\tilde{\ct}(A)$ correspond to a special kind of partial ideal: one which is fixed by every unitary rotation in the sense that $I_{uVu^*} = uI_Vu^*$.  Thus Conjecture \ref{cstarconj} follows from:

\begin{conj}
A partial ideal of a unital $C^*$-algebra arises from a total ideal if and only if it is fixed by every unitary rotation.
\end{conj}

In forthcoming work with Rui Soares Barbosa, we establish the von Neumann algebraic analogue of this conjecture.  One may define a \emph{partial ideal} for a unital von Neumann algebra by replacing  ``unital, commutative sub-$C^*$-algebra''with  ``unital, commutative sub-von-Neumann-algebra'' and ``closed ideal'' with ``ultraweakly closed ideal''.  In this case, a total ideal is an ultraweakly closed, two-sided ideal.

\begin{thm}
A partial ideal of a unital von Neumann algebra arises from a total ideal if and only if it is fixed by every unitary rotation.
\end{thm}

This theorem provides some measure of evidence for the original conjecture's verity and its proof may be adapted to show that the original conjecture holds for a large class of---or perhaps all---$C^*$-algebras.

\section{Acknowledgements}

It is a pleasure to thank my supervisors, Samson Abramsky and Bob Coecke, as well as Chris Heunen, Andreas D\"{o}ring, and Jamie Vicary for their guidance and encouragement during this project.  I also wish to thank Kobi Kremnitzer, George Elliott,  Rui Soares Barbosa, and Brent Pym for many crucial and constructive discussions.  I would like to acknowledge the  support of Merton College, the Oxford University Computing Laboratory, the Clarendon Fund, and NSERC.

\end{document}